\documentclass[11pt]{amsart}

\usepackage[T1]{fontenc}
\usepackage[latin1]{inputenc}
\usepackage{hyperref}

\newcommand\mR{{\mathbb{R}}}

\newcommand\mC{{\mathbb{C}}}
\newcommand\mQ{{\mathbb{Q}}}
\newcommand\mP{{\mathbb{P}}}

\newcommand{\trace}{{\operatorname{tr}}}
\newcommand{\ricci}{{\operatorname{Ric}}}
\newcommand{\real}{{\operatorname{Re}\,}}
\newcommand{\Aut}{{\operatorname{Aut}}}
\newcommand{\diag}{{\operatorname{diag}}}

\newcommand{\fracts}[2]{\frac{\textstyle #1}{\textstyle #2}}
\newcommand{\parpar}[1]{\frac{\partial}{\partial #1}}

\newtheorem{theorem}{Theorem}[section]
\newtheorem{lemma}[theorem]{Lemma}
\newtheorem{corollary}[theorem]{Corollary}

\newtheorem{remark}[theorem]{Remark} 
\newtheorem{definition}[theorem]{Definition} 
\newtheorem{conjecture}[theorem]{Conjecture} 

\newtheorem{example}[theorem]{Example}

\newtheorem{definitionandtheorem}[theorem]{Definition and Theorem}

\begin{document}

\title{Stability and Futaki Invariants of Fano Hypersurfaces} 

\date{November 9, 2001}

\thanks{Supported by Schwerpunktprogramm ``Globale Methoden in der 
        komplexen Geometrie'' Deutsche Forschungsgemeinschaft DFG}

\author[Thomas R. Bauer]{Thomas Rudolf Bauer}

\address{Mathematisches Institut, Universität Bayreuth,
  95440~Bayreuth, Germany} 
\email{thomas.bauer@uni-bayreuth.de}

\maketitle
%\tableofcontents

\begin{abstract}
Let X be a Fano manifold. G.~Tian proves that if X admits a Kähler-Einstein
metric, then it satisfies two different stability conditions: one involving
the Futaki invariant of a special degeneration of X, the other 
Hilbert-Mumford-stability of X w.~r~.t.\ a certain polarization. He conjectures that 
each of these conditions is also sufficient for the existence of such
a metric. If this is true, then in particular the two stability
conditions would be equivalent. We show that for Fano hypersurfaces in
projective space, where due to the work of Lu and Yotov an explicit 
formula for the Futaki invariant is known, these two
conditions are indeed very closely related.
\end{abstract}

\section{Introduction}

For complex manifolds with ample or trivial canonical bundle the existence
of a Kähler-Einstein metric is well known. For Fano manifolds, the situation
is much more complicated. Of course, the Fubini-Study metric on $\mP^n$ 
is Kähler-Einstein. But already $\mP^2$ blown up in 1 or 2 points does 
not admit such a metric. In general, it is supposed that the existence of a 
Kähler-Einstein metric is related to the geometry of the Fano manifold, and 
to be more specific, to certain notions of stability of the manifold.

In \cite{Ti97}, Tian proves that a Fano manifold admitting a Kähler-
Einstein metric satisfies two different stability conditions, one 
involving the Futaki invariant of a special degeneration of the manifold, 
the other Hilbert-Mumford stability of the manifold with respect to a
certain polarization. Tian 
conjectures that each of these conditions is also sufficient for the 
existence of such a metric. If this is true, then in particular the two
stability conditions would be equivalent. We show that for 
Fano hypersurfaces in projective space, where due to the work of Lu 
\cite{Lu99} and Yotov \cite{Yo99} an explicit formula
for the Futaki invariant is known, these two conditions are indeed
closely related, as it turns out that the Futaki invariant $F$ equals Mumford's
$\mu$-function up to a constant.

Furthermore, we introduce the notion of a special degeneration of a
hypersurface \emph{as a hypersurface}, which is very similar to Tian's definition.
Our main theorem is the following

\begin{theorem}
Let $X_f$ be a hypersurface of degree $d$, $1<d<n+1$ in $\mP^n$.
Then the following conditions are equivalent
\begin{itemize}
\item $X_f$ is weakly Hilbert-Mumford stable
\item for every special degeneration $\pi:Y\to\Delta$ of $X_f$ 
as a hypersurface, where the degeneration is induced by a vector field $v$
on $\mP^n$, one has $F_{Y_0}(v_{|Y_0}) \geq 0$ and ``$=0$'' iff the degeneration
is trivial.
\end{itemize}
\end{theorem}

I would like to thank Thomas Peternell for his suggestion to have a closer 
look at Tian's paper \cite{Ti97}, and Thomas Eckl for some helpful comments 
on vector fields and automorphisms.

\section{The Futaki invariant and special degenerations}

We always denote by $X$ a Fano manifold and ask whether 
there exists a Kähler-Einstein metric on $X$, i.~e. a Kähler metric 
such that the associated Kähler form $\omega$ satisfies $\ricci(\omega)=\omega$.
We further denote by $\eta(X)$ the Lie algebra of holomorphic vector fields on $X$.

In \cite{Fu83}, Futaki introduced the character
$$
F_X:\eta(X)\to\mC \qquad\qquad\qquad v\mapsto\int_X v(g)\omega^n
$$
where $n$ is the complex dimension of $X$ and 
$g:X\to\mR$ is a function with $\ricci(\omega)-\omega=i/2\pi \, 
\partial \bar{\partial}g$. He proved that this is independent of 
the choice of the Kähler form $\omega$ with $[\omega]=[c_1(X)]$.
In particular, if $X$ admits a Kähler-Einstein metric then $F_X$
vanishes and therefore the Futaki invariant gives a first obstruction 
to the existence of Kähler-Einstein metrics.

Unfortunately, there exist Fano manifolds without any global holomorphic 
vector fields (i.~e. $F_X\equiv0$ trivially) which do not admit any
Kähler-Einstein metric \cite{Ti97}. Therefore, if we want an equivalent 
condition, we have to refine this.

\begin{definition}\cite{Ti97}
a) A fibration $\pi:Y\to\Delta$ over the unit disc $\Delta$ 
is a \emph{special degeneration of }$X$ 
if 
\begin{itemize}
   \item $\pi$ is smooth over $\Delta-{0}$
   \item $X$ is isomorphic to a fiber $Y_z$ for some $z\in\Delta-{0}$
   \item the special fiber $Y_0$ is a normal variety (in particular irreducible 
and reduced)
   \item the relative anticanonical bundle $-K_{Y/\Delta}$ is ample and 
therefore induces an embedding $Y\subset\mathbb{P}^N\times\Delta$ such that
$\pi$ is induced by the projection on $\Delta$ and
   \item there exists a vector field $w$ on $Y$ with $\pi_*w=-z\parpar{z}$ 
on $\Delta$.
\end{itemize}
b) A special degeneration is said to be \emph{trivial} if $Y=X\times\Delta$,
$\pi$ is the projection on $\Delta$ and $w$ is induced by a vector field on $X$.
\end{definition}

\begin{remark}
1) The last condition of a) implies that $w$ induces a vector field on $Y_0$.
2) In a special degeneration, the special fiber $Y_0$ is always $\mQ$-Fano.
Ding and Tian \cite{DT92} show that the above definition of the Futaki 
invariant still makes sense in this case.
\end{remark}

\begin{theorem}\cite{Ti97}
If $X$ is Kähler-Einstein then for every special degeneration of $X$ one has 
$\real F_{Y_0} (\omega_{|Y_0}) \geq 0$ and ``$=0$'' iff the degeneration is trivial.
\end{theorem}

\begin{conjecture}\cite{Ti97}\label{conj1}
This is also sufficient for the existence of a Kähler Einstein metric on $X$.
\end{conjecture}

Of course, the main difficulty here is the calculation of the Futaki invariant.
If $Y_0$ is smooth or has only orbifold singularities, one can use a fixed 
point formula for $F_{Y_0}$ which was proved by Futaki \cite{Fu88} 
resp.\ Ding and Tian \cite{DT92}. 

If $Y_0$ is a hypersurface in projective space or more general a complete 
intersection, Lu \cite{Lu99} and Yotov \cite{Yo99} independently 
obtained an explicit formula for the Futaki invariant. To keep things easy, 
we will stick to the hypersurface case, although everything what follows 
is very similar in the complete intersection case.

Let $X=X_f=\{f=0\} \subset \mP^n$ be a hypersurface of degree $d$, $1<d<n+1$, and 
denote by $[z_0:\ldots:z_n]$ the coordinates on $\mP^n$. Furthermore, let 
$v=\sum a_{ij} z_j \parpar{z_i}$ be a vector field on $\mP^n$, which is 
normalized by the condition $\sum a_{ii}=0$, and let 
$\{\sigma_t\}_{-\infty<t<\infty}$
be the real 1-parameter subgroup of $\Aut(X)$ generated by $v$.
If $X_f$ is invariant under $\{\sigma_t\}$, then $v(f)=\kappa \cdot f$ with 
$\kappa \in \mC$ and $v$ induces a vector field on $X_f$. By a theorem of 
Bott, every vector field on $X_f$ arises in this way.

\begin{definition}\cite{Lu99}\cite{Yo99}
In the situation above,
$$
\tilde{F}_X(v)=-(n+1-d)(d-1)\fracts{(n+1)}{n}\kappa
$$
is the \emph{generalized Futaki invariant of $v$ on $X$}.
\end{definition}

We remark that $(n+1-d)(d-1)\fracts{(n+1)}{n}$ is positive, as $1<d<n+1$.

\begin{theorem}\cite{Lu99}\cite{Yo99} 
If $X_f$ is $\mQ$-Fano, then $F_X=\tilde{F}_X$.
\end{theorem}

\section{The method of Ding and Tian}

To make things even more explicit, we now recall an example of Ding and Tian, 
illustrating their method of evaluating the Futaki invariant on a degeneration 
of a Fano manifold $X$ to show that $X$ does not admit any Kähler-Einstein 
metric. Again we restrict ourselves to the case of Fano hypersurfaces in $\mP^n$.

Let $X=X_f$ be a hypersurface of degree $1<d<n+1$ in $\mP^n$. If $X$ is smooth and 
$d \neq 2$ then $\eta(X)=0$ by a theorem of Kodaira and therefore $F_X \equiv 0$ 
trivially. Let $v$ be a vector field on $\mP^n$ and let $\{\sigma_t\}_{-\infty<t<\infty}$ 
be the induced real 1-parameter subgroup of $\Aut(\mP^n)$. Furthermore, let 
$X_t=\sigma_t(X)$ and let $X_\infty=\lim_{t\to\infty}X_t$ be the limit of the 
$X_t$, if such a limit exists. We also assume that $X_\infty$ is a normal variety. 
By construction, $v_{|X_\infty}$ is a vector field on $X_\infty$. The main theorem 
of \cite{DT92} now asserts that if $X$ is Kähler-Einstein, then one has 
$\real F_{X_\infty}(v_{|X_\infty})\geq0$, where on $X_\infty$ one uses again the 
generalized Futaki invariant. This result holds in the orbifold case as well 
(by using local uniformization).

\begin{example}
Let $f=z_0 z_1^2 + z_2 z_3 (z_2-z_3) + z_1 f_2(z_1,z_2,z_3)$ where $f_2$ is 
a homogeneous polynomial of degree 2 in $z_1,z_2,z_3$, and $X=X_f=\{f=0\}$ the 
corresponding cubic surface in $\mP^3$ which has a $D_4$-singularity.
Let $v=-7z_0\parpar{z_0}+5z_1\parpar{z_1}+z_2\parpar{z_2}+z_3\parpar{z_3}$ which 
already has trace $0$. One easily reads off that $\sigma_t.z_0=e^{-7t}z_0$, 
 $\sigma_t.z_1=e^{5t}z_1$, $\sigma_t.z_2=e^tz_2$ and $\sigma_t.z_3=e^tz_3$ so 
that $\sigma_t.z_0 z_1^2 = e^{-3t}z_0 z_1^2$ and $\sigma_t.z_2 z_3 (z_2-z_3) = 
e^{-3t}z_2 z_3 (z_2-z_3)$ whereas for any other monomial $m$ in $f$ we get 
$\sigma_t.m=e^{-at}m$ with $a>3$. Therefore the limit is $X_\infty = 
\{z_0 z_1^2 + z_2 z_3 (z_2-z_3) =0\}$ which is still a Fano orbifold.
Moreover, $v(f_\infty)=3f_\infty$, where $X_\infty$ is given by $f_\infty$, and 
with the formula of Lu and Yotov we calculate that $F_{X_\infty}(v_{|X_\infty}) = 
-1\cdot2\cdot\frac{4}{3}\cdot3=-8<0$. Therefore $X_f$ does not admit a Kähler-Einstein 
(orbifold) metric.
\end{example}

We remark that $X_f$ is Hilbert-Mumford unstable, which of course is already 
in \cite{GIT}, but actually the calculation is the same what we did above. 
We will come back to this later.

With some more work, Ding and Tian actually prove that if a normal cubic in $\mP^3$
admits a Kähler-Einstein orbifold metric then it is Hilbert-Mumford semistable.

\section{Tian's 2nd stability condition}

Now we recall Tian's second stability condition which is necessary for a Fano 
manifold to allow a Kähler-Einstein metric. If we state the result for the 
hypersurface case only, it is almost trivial, as any nonsingular 
hypersurface of degree $>2$ in $\mP^n$ is Hilbert-Mumford stable \cite{GIT}.
Nevertheless we belive that it is more pedagogic to ignore this fact
and state only a trivial case of a highly nontrivial theorem.

\begin{theorem}\cite{Ti97}
Let $X_f$ be a smooth hypersurface of degree $d$ in $\mP^n$, $1<d<n+1$. If $X_f$ is 
Kähler-Einstein, then $X_f$ is weakly Hilbert-Mumford stable, i.~e. the orbit of 
$f\in\mC^N$, $N=\binom{n+d}{d}$ under the natural $SL(n+1)$-action is closed. If in 
addition $\eta(X_f)=0$ then $X_f$ is Hilbert-Mumford stable.
\end{theorem}

\begin{conjecture}\cite{Ti97}
This is also sufficient for the existence of a Kähler-Einstein metric.
\end{conjecture}

In the general case, Tian proves that $X$ is weakly stable w.~r.~t. the polarization 
given by a certain bundle, cf.\ \cite{Ti97}.

\section{Stability of hypersurfaces}

Once more, we fix the notation for the following statements. Let $v$ be a vector field 
on $\mP^n$. Via the natural projection $v$ is induced by a vector field 
$\sum_{i,j=0}^n a_{ij} z_j \parpar{z_i}$ on $\mC^{n+1}$. We will not distinguish 
between these two vector fields below. Let $A$ be the matrix $(a_{ij})$ where we 
assume that $v$ is normalized by the condition $\trace(A)=0$. After a linear change 
of variables, we may assume that $A$ is in Jordan form with blocks
$$
\begin{pmatrix}
\lambda_i & 1 & & \\
& \ddots & \ddots & \\
& & \lambda_i & 1 \\
& & & \lambda_i
\end{pmatrix}
$$
We call $v$ a \emph{real} vector field if $A$ has only real eigenvalues, i.~e.\ 
$\lambda_i\in\mR$ for every Jordan block of $A$. Of course a real vector field is still
holomorphic. We note that for real vector fields a limit $X_\infty$ in the sense of Ding 
and Tian always exists.

If the first Jordan block of $A$ operates on the coordinates $z_0,\ldots,z_r$ then the 
induced 1-parameter subgroup $\{\sigma_t\}$ of $\Aut(\mC^{n+1})$ acts on 
$(z_0,\ldots,z_r)^t$ by multiplication with the matrix
$$
\begin{pmatrix}
e^{\lambda_1 t} & t e^{\lambda_1 t} & \frac{t^2}{2} e^{\lambda_1 t} & \cdots 
& \frac{t^r}{r!} e^{\lambda_1 t} \\
 & e^{\lambda_1 t} & t e^{\lambda_1 t} & & \vdots \\
 & & \ddots & \ddots & \vdots \\
 & & & e^{\lambda_1 t} & t e^{\lambda_1 t} \\
 & & & & e^{\lambda_1 t}
\end{pmatrix}
$$
(if we regard the $z_i$ as hyperplanes, then $\{\sigma_t\}$ acts by multiplication
with the inverse matrix as usual). In particular, if $A$ is a diagonal matrix with 
entries $(\lambda_0,\ldots,\lambda_n)$ we get a diagonal action $\sigma_t . z_i = 
e^{\lambda_i t} z_i$. If all the $\lambda_i$ are integers, this induces a $\mC^*$-
action on $\mC^{n+1}$ resp.\ $\mP^n$. By the explicit description above, we also 
note that $\trace(A)=0$ corresponds to $\det(\sigma_t)=1$, i.~e.\ $\sigma_t\in 
SL(n+1)$. 

Finally, $X_f=\{f=0\}$ is the hypersurface in $\mP^n$ given by the homogeneous polynomial 
$f=\sum_{|\gamma|=d} f_\gamma z^\gamma$ of degree $d$, $1<d<n+1$, with 
$\gamma=(\gamma_0,\ldots,\gamma_n)$ a multiindex. We also denote by $f$ the vector 
$(f_\gamma)$ in $\mC^N$.

We recall now some of the definitions and theorems of Mumford's \cite{GIT}. As many 
assertions are quite simple in the hypersurface case, we could not resist to give 
some of the proofs here. In some way this shows how everything fits together very 
nicely.

From Mumford we already know the $\mu$-function, which makes sense for any real 
diagonal vector field $v$ as above with entries $\lambda$ .

\begin{definition}
$\mu_v (f) = \min \{ \lambda \cdot \gamma \, | \, f_\gamma \neq 0 \}$
\end{definition}

If we call $\lambda \cdot \gamma$ the weight of $z^\gamma$ then $\mu_v(f)$ 
is the least weight in $f$.

\begin{definitionandtheorem}\cite{GIT}
a) We call $f$ \emph{Hilbert-Mumford-stable} if the orbit of $f\in\mC^N$ 
under the natural $SL(n+1)$-action is closed and the stabilizer of $f$ is finite.
Equivalent conditions are
\begin{itemize}
\item for all algebraic 1-parameter subgroups of $SL(n+1)$ $f$ has both positive 
and negative weights
\item for all $\mC^*$-vector fields $v\neq0$ on $\mC^{n+1}$ with trace $0$ one has 
$\mu_v(f)<0$.
\end{itemize}
b) We say that $f$ is \emph{weakly Hilbert-Mumford stable} if the orbit of 
$f$ is closed. Equivalently, for all $\mC^*$-vector fields $v\neq0$ with trace $0$ either 
$\mu_v(f)<0$ or the $\mC^*$-orbit of $f$ is a point, i.~e.\ $\mu_v(f)\leq0$ and ``$=0$'' iff 
$f$ is a fixed point of the action induced by $v$. 
\end{definitionandtheorem}

\begin{lemma}
Let $v=\sum \lambda_i z_i \parpar{z_i}$ be a diagonal vector field and assume that 
$v(f)=\kappa f$ for some $\kappa\in\mC$. Then $\mu_v(f)=\kappa$.
\end{lemma}

\begin{proof}
As $v(z^\gamma) = \sum_{i=0}^n \lambda_i z_i \gamma_i z^{\gamma-e_i} = (\lambda \cdot 
\gamma) z^\gamma$, $e_i$ the $i$-th unit vector, we have $\kappa f = v(f) = 
\sum_\gamma (\lambda \cdot \gamma) f_\gamma z^\gamma$ and thus $\lambda \cdot \gamma = 
\kappa$ for all $\gamma$ where $f_\gamma \neq 0$, in particular $\mu_v(f)=\kappa$.
\end{proof}

\begin{lemma}
Assume that $v$ is a diagonal real vector field and let $\mu=\mu_v(f)$.
Then $X_\infty=\{f_\infty=0\}$ is defined by $f_\infty=\sum_{\lambda\cdot\gamma=\mu} 
f_\gamma z^\gamma$ and therefore $v(f_\infty)=\mu f_\infty$ with $\mu=\mu_v(f)= 
\mu_v(f_\infty)$.
\end{lemma}

\begin{proof}
For any $\gamma$ with $f_\gamma \neq 0$ one has 
$\lambda \cdot \gamma - \mu \geq 0$ by the definition of $\mu$.
In the limit, all monomials $z^\gamma$ where $\lambda \cdot \gamma - \mu > 0$ vanish, 
i.~e.\ we only keep those monomials where the minimum $\mu$ is obtained.
\end{proof}

\begin{corollary}
Let $X_f$ be a hypersurface of degree $1<d<n+1$. Then $X_f$ is Hilbert-Mumford stable 
iff for every $\mC^*$-vector field $v\neq0$ with trace $0$ one has $F_{X_\infty}(v)>0$.
\end{corollary}

We emphasize that we make no conditions on $X_\infty$ here.

\begin{proof}
We just note that by the formula of Lu and Yotov,
$$
F_{X_\infty}(v)=-(n+1-d)(d-1)\fracts{n+1}{n}\kappa
$$ where $(n+1-d)(d-1)\fracts{n+1}{n}>0$ and $\kappa=\mu_v(f)$.
\end{proof}

\begin{lemma}
$X_f$ is Hilbert-Mumford stable iff for every diagonal real vector field $v\neq0$ with 
trace $0$ one has $\mu_v(f)<0$.
\end{lemma}

\begin{proof}
We only have to proof the ``only if'' direction. Assume that there exists a real vector 
field $0\neq v = \sum \lambda_i z_i \parpar{z_i}$ with $\sum \lambda_i =0$ and $\mu_v(f)=
\min\{ \lambda \cdot \gamma \, | \, f_\gamma \neq 0 \} \geq 0$, i.~e.\ the 
following system of linear inequalities has a real solution under the constraint 
$\sum \lambda_i =0$: 
$$
\forall \gamma \text{ with } f_\gamma \neq 0: \qquad \lambda \cdot \gamma \geq 0
$$
As all coefficients $\gamma_i$ are nonnegative integers, this system then also has a 
rational solution. By clearing denominators, we obtain a $\mC^*$-vector field with 
$\mu\geq0$. This is a contradiction to our assumption that $f$ is stable.
\end{proof}

This also shows how to calculate destabilizing vector fields.

If $v$ is a real Jordan vector field, we can decompose it into the diagonal part $v_\diag$ 
and the nilpotent part $v^+$: $v=v_\diag + v^+$. The diagonal part induces a $(\mC^*)^k$-
action for some $k$, whereas the nilpotent part induces a $\mC^+$-action. If $X_f$ is 
invariant under $v$ then also under $v_\diag$ and $v^+$.

\begin{lemma}
If $v$ is nilpotent and $X_f$ invariant under $v$ then $v(f)=0$ (i.~e.\ $\kappa=0$).
\end{lemma}

\begin{proof}
As $X_f$ is invariant under $v$, $v(f)=\kappa f$ for some $\kappa$. As $v$ is nilpotent, 
we know that $v=z_1\parpar{z_0} + \cdots + z_r\parpar{z_{r-1}}+\cdots$ for some $r$, where 
we have written out the terms corresponding to the first Jordan block. If $z^\gamma$ is 
the least monomial occuring in $f$ for the lexicographical order, we note that $v(f)$
does not contain $z^\gamma$ and therefore $\kappa=0$.
\end{proof}

\begin{remark}
This is only a very special case of a theorem of Mabuchi \cite{Ma90}:
The Futaki invariant vanishes on nilpotent vector fields.
\end{remark}

So we see that if we degenerate $X_f$ with a $\mC^+$-vector field $v$ we only get 
$F_{X_\infty}(v_{|X_\infty})=0$. But if $X_f$ is Kähler-Einstein we wanted 
to prove some kind of stability and therefore expected to get ``$>0$''. As a result, 
we restrict ourselves to diagonal vector fields (more exactly to vector fields $v$ with 
nonvanishing diagonal part $v_\diag$, but as $F(v^+)=0$ the diagonal part is 
essential).

We sum this up in the following

\begin{corollary}
Let $X_f$ be a hypersurface in $\mP^n$ of degree $1<d<n+1$. The following conditions are 
equivalent:
\begin{itemize}
\item $X_f$ is Hilbert-Mumford stable
\item for all real Jordan vector fields $v$ with trace $0$ one has 
$F_{X_\infty}(v_{|X_\infty}) \geq 0$ and ``$=0$'' only if $v_\diag=0$.
\end{itemize}
\end{corollary}

\section{Special degenerations of hypersurfaces}

Again $X_f=\{f=0\}$ is a hypersurface in $\mP^n$ of degree $d$, $1<d<n+1$.

\begin{definition}
a) A \emph{special degeneration of $X_f$ as a hypersurface} is a fibration 
$\pi:Y\to\mC$ such that
\begin{itemize}
\item $Y$ is a hypersurface in $\mP^n \times \mC$ and $\pi$ is the 
restriction of the projection to the second factor
\item for all $s\in\mC$ the fiber $Y_s$ is a hypersurface of degree $d$ 
in $\mP^n \times \{s\} \cong \mP^n$ and $Y_1 = X_f$
\item there exists a vector field $v$ on $\mP^n$ such that $Y$ is invariant
under $v-s\parpar{s}$.
\end{itemize}
b) A special degeneration of $X_f$ as a hypersurface is said to be 
\emph{trivial} if $Y_s$ is invariant under $v$ for one (any) $s\neq0$ and 
therefore $Y = X_f \times \mC \subset \mP^n \times \mC$.
\end{definition}

\begin{remark}
If we compare our definition with Tian's there are two major differences:
First, we assume that the total space of the degeneration is a hypersurface 
in $\mP^n \times \mC$. Secondly, we do \emph{not} assume that the special 
fiber is a normal variety. In our setting the special fiber is an arbitrary 
hypersurface in $\mP^n \times \{0\}$. We don't know whether it is always 
possible to construct destabilizing degenerations in such a way that 
the special fiber is normal, irreducible and reduced.
\end{remark}

\begin{theorem}
A special degeneration of $X_f$ as a hypersurface is uniquely determined 
by $v$.
\end{theorem}

\begin{proof}
We put $w=v-s\parpar{s}$. We assume that $Y$ is given by the polynomial 
$H(s)$ where $H(s_0)$ is a homogeneous polynomial of degree $d$ in 
$z_0,\ldots,z_n$ for each $s_0\in \mC$, in particular $H(1)=f$. As $Y$ is 
invariant under $w$ we have $w(H(s))=\kappa H(s)$ for some constant $\kappa$.

If we put $u=z_0\parpar{z_0}+\ldots+z_n\parpar{z_n}$ then $u(H(s))=d H(s)$ 
by Euler's formula, so if we modify $w$ by a suitable multiple of $u$ we
can assume that $\kappa=0$. We note that hereby we have lost our former 
normalization that $v$ has trace $0$.

Integration of $v$ gives the real 1-parameter subgroup $\{\sigma_t\}$ of 
$\Aut(\mC^{n+1})$. We put $F(t)=\sigma_t . f$. Then for any $t\in\mR$ clearly 
$F(t)$ defines a hypersurface in $\mP^n \times \{t\}$. By construction 
(as $(z\parpar{z} + \parpar{t}) (e^{-t} z) = 0$ and so on) we know that
$(v+\parpar{t})(F(t))=0$, i.~e.\ $\{F(t)=0\} \subset \mP^n \times \mR$ 
is invariant under $v+\parpar{t}$.

Now we make the substitution $s:\mR \to \mR^+$, $s=e^{-t}$ which gives 
$\parpar{t} = -e^{-t}\parpar{s} = -s\parpar{s}$. Then $G(s):=F(-\ln(s))$
is invariant under $v-s\parpar{s}=w$, more precisely $w(G(s))=0$ and 
furthermore $G(1)=F(0)=f=H(1)$.

From the uniqueness theorem for the solution of a differential equation 
we conclude now that $Y=\{G(s)=0\}$ over $\mR^+$ and as everything is 
analytic, therefore over $\mC$.
\end{proof}

From the fact that $G(s)$ extends to a holomorphic function on $\mC$ we
derive the following corollaries.

\begin{corollary}
In a special degeneration as a hypersurface the nilpotent part $v^+$
acts trivially on $f$.
\end{corollary}

\begin{proof}
First we consider only the nilpotent part $v^+$. Let $\{\psi_t\}$ be the 
real 1-parameter subgroup generated by $v^+$. Then we get
$$
F^+ (t) := \psi_t . f = f + t f_1 + t^2 f_2 + \ldots + t^m f_m
$$
for some $m$ where the $f_i$ are polynomials in $z_0,\ldots,z_n$.
Consequently
$$
G^+ (s) = f - (\ln s) f_1 + (\ln s)^2 f_2 - + \ldots
$$
which is only holomorphic on $\mC$ if $f_1=\ldots=f_m=0$, i.~e.\
$\psi_t . f = f$.
In the general case we get additionally some terms of the form 
$e^{at}$ which give some powers of $s$ which don't make the whole 
thing holomorphic unless all of the logarithmic expressions vanish.
\end{proof}

As the nilpotent part acts trivial, we may assume now that $v$ is 
a diagonal vector field.

\begin{corollary}
In a special degeneration as a hypersurface $v$ acts $\mC^*$-like 
on $f\in\mP^{N-1}$.
\end{corollary}

\begin{proof}
As
$$
\sigma_t . z^\gamma = e^{-(\lambda \cdot \gamma) t} z^\gamma = 
s^{\lambda \cdot \gamma} z^\gamma
$$
we calculate that
$$
G(s)=\sum f_\gamma s^{\lambda \cdot \gamma} z^\gamma
$$
which is only holomorphic on $\mC$ if $\lambda \cdot \gamma$ is 
a nonnegative integer for all $\gamma$ with $f_\gamma \neq 0$.
\end{proof}

We recall that we cannot assume that $v$ has trace $0$ here, but 
we normalized $v$ such that $w(H)=0$.

\begin{corollary}
There exists a vector field $v'=\sum \lambda_i' z_i \parpar{z_i}$ 
with $\lambda_i \in \mQ$ which acts on $f$ as $v$ does.
\end{corollary}

\begin{proof}
Complete $\{1\}$ to a basis of $\mC$ as a vector space over $\mQ$ and 
write $\lambda_i = \lambda_i' \cdot 1 + \ldots$ with coefficients 
$\lambda_i' \in \mQ$. We know that all monomials $z^\gamma$ in $f$ 
have integer weight $\lambda \cdot \gamma$ for the action of $v$.
Expressing the weights in terms of our basis $\{1,\ldots\}$ we 
conclude that only $v'=\sum \lambda_i' z_i \parpar{z_i}$ acts 
nontrivially on $f$.
\end{proof}

Of course we cannot assume that $v'$ has integer coefficients, as we 
see by the stupid example of the action induced by $\frac{1}{2}z 
\parpar{z}$ on the monomial $z^2$. If we want to do so in our 
special degeneration as a hypersurface, we must allow a rescaling on 
the $\mC$-factor, i.~e. we have to multiply $-s\parpar s $ by a 
suitable integer.

Without loss of generality we may therefore assume that 
$v=\sum \lambda_i z_i \parpar{z_i}$ with $\lambda_i \in \mQ$.

\begin{theorem}
Let $X_f$ be a hypersurface of degree $d$, $1<d<n+1$ in $\mP^n$. 
Then the following conditions are equivalent
\begin{itemize}
\item $X_f$ is weakly Hilbert-Mumford stable
\item for all special degenerations of $X_f$ as a hypersurface one has 
$F_{Y_0} (v_{|Y_0}) \geq 0$ and ``$=0$'' iff the degeneration is trivial.
\end{itemize}
\end{theorem}

\begin{proof}
Let us first assume that $f=\sum f_\gamma z^\gamma$ is weakly 
Hilbert-Mumford stable and that $\pi:Y\to\mC$ is a special degeneration 
of $X_f$ as a hypersurface determined by $v$. As we noted above we may 
assume that $v$ is a diagonal vector field with rational coefficients. 
Let $v'$ be the normalization of $v$, i.~e.\ we add to $v$ a suitable 
multiple of $u=\sum z_i \parpar{z_i}$ to get trace $0$. Furthermore let
$v''=M v'$ where $M$ is a suitable positive integer such that $v''$ is 
a $\mC^*$-vector field. We note that the Futaki invariants of $v$ and 
$v''$ have the same sign. 

As $f$ is weakly Hilbert-Mumford 
stable we have $\mu_v(f) \leq 0$. Because the special fiber 
of a special degeneration cleary corresponds to the limit variety 
$X_\infty$ in the sense of Ding and  Tian (using our change of variables 
$s=e^{-t}$), we conclude using the results of the previous paragraph 
that $F_{Y_0} ({v''}_{|Y_0}) \geq 0$ and if we already have ``$>0$'' we are 
done. 

So we may further assume that $F_{Y_0}({v''}_{|Y_0}) =0$ and 
therefore $\mu_{v''}(f)=0$. As $f$ is weakly Hilbert-Mumford stable this 
is possible only if $f\in\mC^N$ is a fixed point of the $SL(n+1)$-action,
i.~e.\ if $\{\sigma''_t\}$ is the 1-parameter subgroup generated by $v''$ 
then $\sigma''_t . f = f$. Let $\sigma_t$ be the 1-parameter subgroup 
generated by $v$. As $v$ and $v'$ differ only by a multiple of $u$ we 
conclude that $\sigma_t . f = e^{-at} f$ for some constant $a$. This 
means that $Y$ is defined by the equation $s^a f=0$. As the special 
fiber is a hypersurface in $\mP^n \times \{0\}$ as well, we conclude 
that $a$ must be $0$ and the degeneration is trivial (consequently, 
already $v$ had trace $0$ in this case).

If the degeneration is trivial, we have $\sigma_t . f = f$ and therefore
$\sigma''_t . f = e^{-at} f$. As $f$ is weakly Hilbert-Mumford stable we know
that $a=\mu_{v''}(f) \leq 0$. But as $\sigma''_{-t} . f = e^{a(-t)}$ for 
the action of the 1-parameter subgroup induced by $-v''$ also 
$-a=\mu_{-v''}(f) \leq 0$ and therefore $a=0$.

To prove the converse, we assume that $f$ fulfills the second condition.
Let $v$ be a $\mC^*$-vector field with trace $0$ and $\mu=\mu_v(f)$. 
Let $v'=d v - \mu u$
and $\{\sigma_t\}$, $\{\sigma'_t\}$ the 1-parameter subgroups induced by
$v$ resp.\ $v'$. If we sort the terms in $\sigma_t . f$ according to 
their weight we get
$$
\sigma_t . f = e^{-\mu t} f_0 + e^{-(\mu + 1) t} f_1 + 
e^{-(\mu + 2) t} f_2 +\ldots
$$ for some homogeneous polynomials $f_i$ of degree $d$, and as $u(f)=d f$
$$
\sigma'_t . f = f_0 + e^{-dt} f_1 + e^{-2dt} f_2 + \ldots
$$
and $v'$ induces a special degeneration of $X_f$ as a hypersurface. By 
our assumption $F_{Y_0}({v'}_{|Y_0}) \geq 0$ and therefore $\mu_v(f) \leq 0$.

If $\mu_v(f)=0$ even $v$ induces a special degeneration which must be 
trivial. This means that $\sigma_t . f = f$ and we are done. If $f$ is 
fixed under the operation of $\mC^*$ induced by $v$, again $v$ induces 
the trivial special degeneration and therefore $\mu_v(f)=0$
\end{proof}

\begin{example}
Let $f=z_0 z_1^2 + z_2 z_3 (z_2 - z_3) + z_1 f_2 (z_1,z_2,z_3)$. Then 
$v=-7 z_0 \parpar{z_0} +5 z_1 \parpar{z_1} + z_2 \parpar{z_2} + 
z_3 \parpar{z_3}$ is a destabilizing vector field with $\mu_v(f)=3$. If we 
define $v'=-24 z_0 \parpar{z_0} +12 z_1 \parpar{z_1}$ then $v'$ induces a
special degeneration of $X_f$.
\end{example}


\begin{thebibliography}{DT92}

\bibitem[DT92]{DT92} W.~Ding, G.~Tian: \emph{Kähler-Einstein metrics and the 
generalized Futaki invariant.} Inv. Math. 110, 315--335 (1992)

\bibitem[Fu83]{Fu83} A.~Futaki: \emph{An obstruction to the existence of 
Einstein Kähler metrics.} Inv. Math. 73, 437--443 (1983)

\bibitem[Fu88]{Fu88} A.~Futaki: \emph{Kähler-Einstein metrics and integral 
invariants.} Springer LNM 1314 (1988)

\bibitem[Lu99]{Lu99} Z.~Lu: \emph{On the Futaki invariants of complete
intersections.} Duke Math. J. 100, 359--372 (1999)

\bibitem[Ma90]{Ma90} T.~Mabuchi: \emph{An algebraic character associated with 
the Poisson brackets.} Adv. Studies Pure Math. 18-I, 339--358 (1990)

\bibitem[GIT]{GIT} D.~Mumford, J.~Fogarty: \emph{Geometric invariant theory.
2nd enlarged edition.} Springer-Verlag Berlin Heidelberg New York (1982) 

\bibitem[Ti94]{Ti94} G.~Tian: \emph{The K-energy of hypersurfaces and 
stability.} Commun. Anal. Geom. 2, 239--265 (1994)

\bibitem[Ti97]{Ti97} G.~Tian: \emph{Kähler-Einstein metrics with positive scalar 
curvature.} Inv. Math. 130, 1--37 (1997)

\bibitem[Ti00]{Ti00} G.~Tian: \emph{Canonical metrics in Kähler geometry.}
Birkhäuser Verlag Basel, Berlin, Boston (2000)

\bibitem[Yo99]{Yo99} M.~Yotov: \emph{On the generalized Futaki invariant.}
LANL-preprint math.AG/9907055 (July 1999)
  
\end{thebibliography}
\end{document}